\documentclass[journal]{IEEEtran}
\usepackage[utf8]{inputenc}
\usepackage{cite}
\usepackage{graphicx}
\usepackage{url}
\usepackage{lipsum}
\usepackage{color}
\usepackage{xcolor} 
\usepackage{amsmath}
\usepackage[english]{babel}
\usepackage{amsthm}

\theoremstyle{definition}

\usepackage{mathtools}
\usepackage{array}
\usepackage{multirow}
\usepackage{makecell}
\usepackage{float}
\usepackage{setspace}
\usepackage{supertabular,booktabs}
\usepackage{amssymb}
\usepackage{algorithm}
\usepackage{physics}
\usepackage[shortlabels]{enumitem}
\usepackage{algpseudocode}
\usepackage{comment}
\newboolean{showcomments}
\setboolean{showcomments}{true}

\usepackage{amsfonts}
\usepackage{multicol}
\usepackage[top=1.2 cm, bottom=1.2 cm, left=1.35 cm, right=1.35 cm]{geometry}

\allowdisplaybreaks[4]
\usepackage{authblk}
\usepackage{nomencl}
\usepackage[symbol]{footmisc}

\usepackage{subcaption}
\newcommand{\changed}[1]{\textcolor{black}{#1}}

\newcommand\figWtwo{0.241}

\newcommand{\julio}[1]{\ifthenelse{\boolean{showcomments}}
        { \textcolor{red}{(JB:  #1)}}{}}
\newcommand{\brian}[1]{\ifthenelse{\boolean{showcomments}}
        { \textcolor{blue}{(BL:  #1)}}{}}
\makenomenclature
\usepackage[colorlinks=true,linkcolor=black,anchorcolor=black,citecolor=black,urlcolor=black]{hyperref}

\begin{document}

\title{\LARGE Robust Dynamic Operating Envelopes via Superellipsoid-based Convex Optimisation in Unbalanced Distribution Networks}

\author{Bin Liu,~\IEEEmembership{Member,~IEEE} and Julio H. Braslavsky,~\IEEEmembership{Senior Member,~IEEE}
	\thanks{Bin Liu (\textit{corresponding author}) was with Energy Centre, Commonwealth Scientific and Industrial Research Organisation (CSIRO), Newcastle 2300, Australia. He is now with the Network Planning Division, Transgrid, Sydney 2000, Australia. Julio, H. Braslavsky is with Energy Centre, CSIRO, Newcastle 2300, Australia. (Email: eeliubin@ieee.org, julio.braslavsky@csiro.au)}}

\markboth{Journal of \LaTeX\ Class Files,~Vol.~17, No.~8, December~2022}%
{Shell \MakeLowercase{\textit{et al.}}: A Sample Article Using IEEEtran.cls for IEEE Journals}

\maketitle

\begin{abstract}
Dynamic operating envelopes (DOEs) have been introduced to integrate distributed energy resources (DER) in distribution networks via real-time management of network capacity limits. Recent research demonstrates that uncertainties in DOE calculations should be carefully considered to ensure network integrity while minimising curtailment of consumer DERs. This letter proposes a novel approach to calculating DOEs that is robust against uncertainties in the utilisation of allocated capacity limits and demonstrates that the reported solution can attain close to global optimality performance compared with existing approaches.
\end{abstract}

\begin{IEEEkeywords}
	DER integration, dynamic operational envelopes, superellipsoid, unbalanced optimal power flow, hosting capacity.
\end{IEEEkeywords}

\section{Introduction}
Dynamic operating envelopes (DOEs) specify the operational capacity range for consumer DERs that is permissible at their connection point and appear as a key enabler of emerging flexible power system architectures and have gained growing interest within industry and academia as an instrument to efficiently manage DER export/import limits \cite{DEIP2022,Liu2021-doe}. With substantial advances being made in developing approaches to calculating DOEs, recent publications revealed that uncertainties arising from parameter or forecast accuracy and true utilisation of DOEs need to be carefully catered for to make the calculated DOEs more robust \cite{liu2022robust,BL_ieee_access,Yi2022}.
\begin{figure}[htbp!]
	\centering
    	\begin{subfigure}[b]{\figWtwo\textwidth}
     		\centering\includegraphics[width=\textwidth]{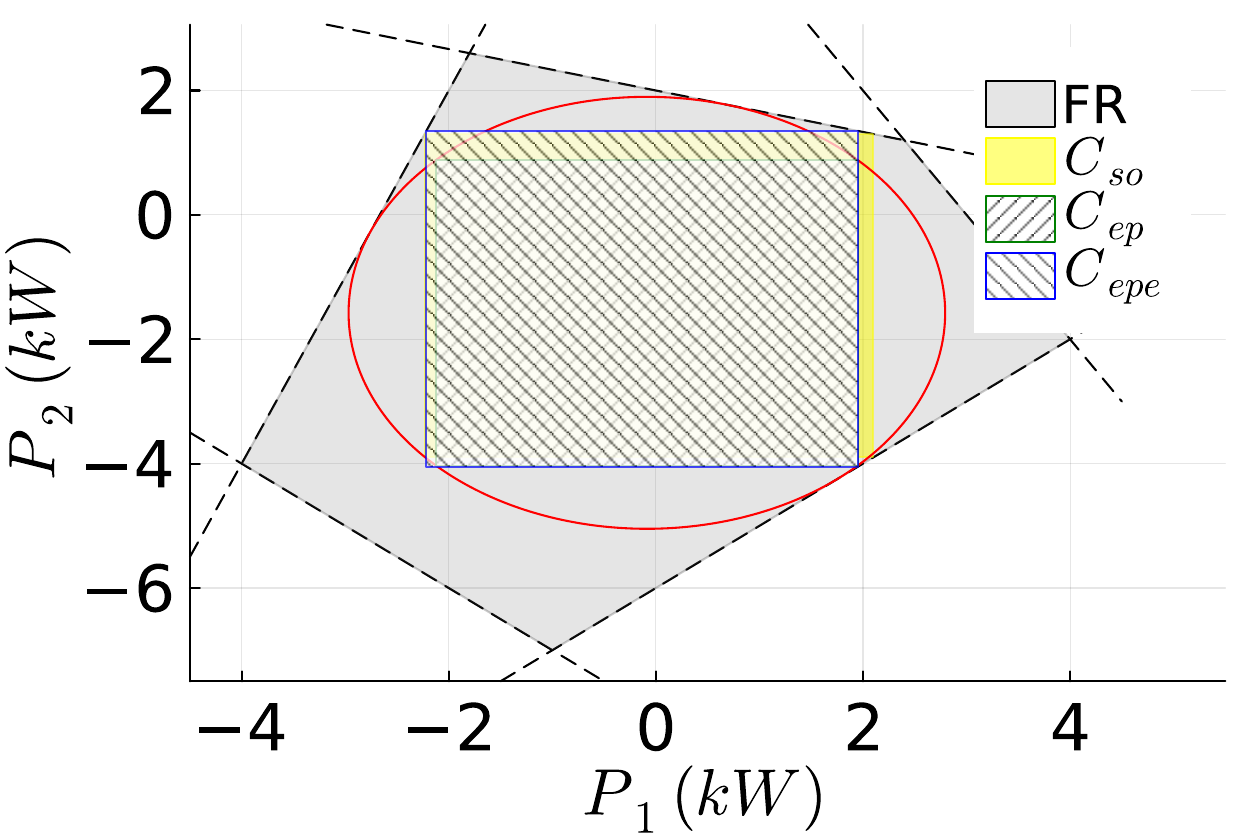}
    	\caption{\small RDOE concept.}
         \label{fig_RDOE_concept}
     \end{subfigure}
	\hfill
    	\begin{subfigure}[b]{\figWtwo\textwidth}
     		\centering\includegraphics[width=\textwidth]{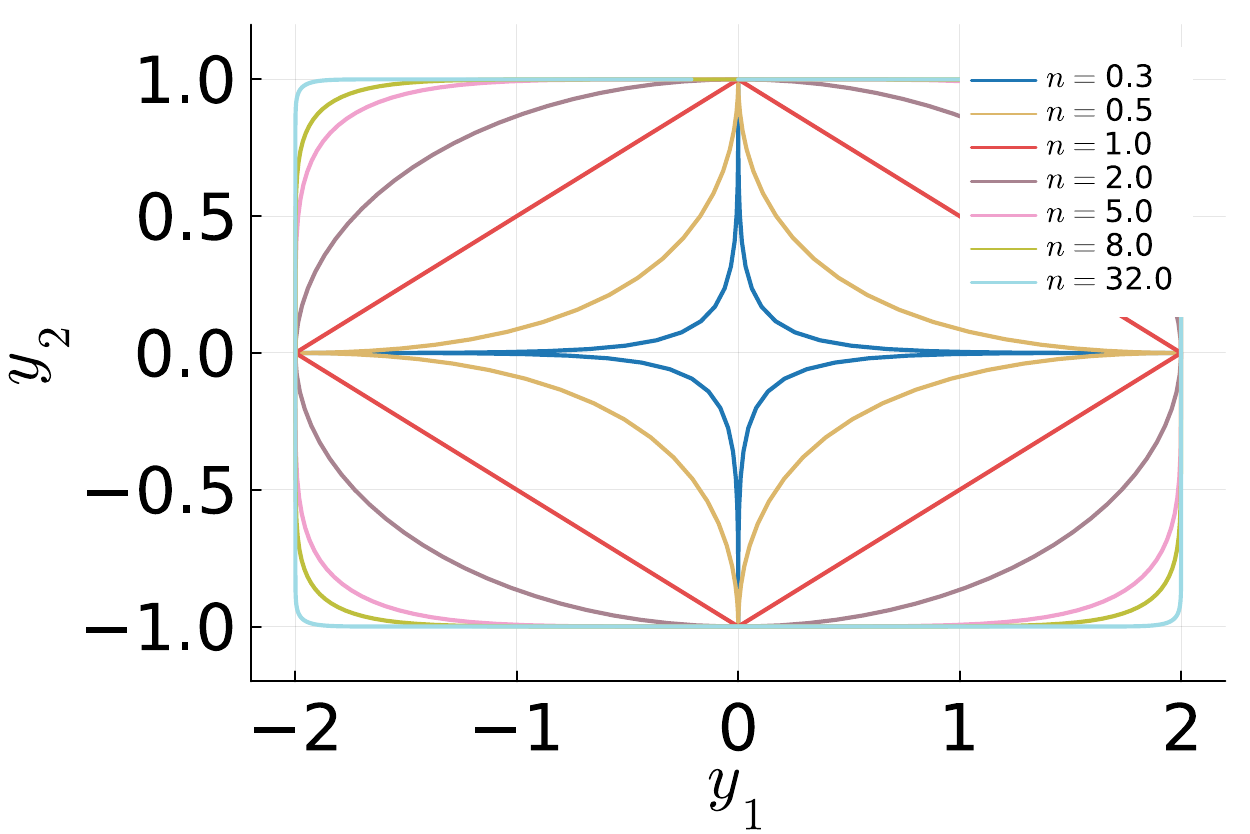}
    	\caption{\small Superellipsoid concept.}
              \label{fig_superellipsoid}
     \end{subfigure}
	\hfill
	\caption{\small (a) Conceptual example in calculating RDOEs via three steps; (b) Superellipsoid with $u_c=[0,0]^T$ and $L=[2,0;0,1]$.}
    \label{fig_RDOE_concept_hypeellipse}
\end{figure}

In \cite{liu2022robust}, a linear unbalanced three-phase optimal power flow (UTOPF)-based approach was proposed to calculate the robust DOEs (RDOEs), which can explicitly hedge against uncertainties in active customers' utilisation of their allocated capacities. 
The RDOE concept can be explained by the feasible region (FR) for DERs \cite{Riaz2022}, as shown in Fig.\ref{fig_RDOE_concept}. For the two active customers, even if the allocated DOEs are $(-4~\text{kW},~-4~\text{kW})$, operational violations can still occur, noting that any operational point $(P_1,P_2)$ with $P_1=-4~\text{kW}$ and $P_2>-4~\text{kW}$ will be out of the FR. 
To address the robustness issue, a three-step approach capable of considering the operational statuses of active customers and optimising controllable reactive powers was proposed in \cite{liu2022robust}. As shown in Fig.\ref{fig_RDOE_concept}, the approach includes: 1) seeking the maximum hyperrectangle ($\mathcal{C}_{ep}$), measured as its volume, within the FR as a strategy to \emph{decouple} the interdependence of DOEs for all active customers via the maximum inscribed hyperellipsoid (\emph{red} curve), which also reports optimal dispatched reactive powers; 2) With fixed optimal reactive dispatch, removing redundant constraints from the updated FR; and 3) Expanding the $\mathcal{C}_{ep}$ based on Motzkin Transposition Theorem (MTT), leading to $\mathcal{C}_{epe}$ and the RDOEs.

\color{black}
However, compared with the stochastic optimisation (SO)-based approach ($\mathcal{C}_{so}$ in Fig.\ref{fig_RDOE_concept}), sub-optimality in the three-step approach is inevitable due to optimising reactive powers in the first step only and the non-convex formulation in the final step. To further improve the optimality of the solution, this letter proposes a novel approach that can report a solution close enough to \emph{global} optimality. The contributions of this letter are:
1) A convex formulation is proposed based on a \emph{superellipsoid} (instead of a \emph{hyperellipsoid} in the first step in the three-step approach) to seek the \emph{maximum} hyperrectangle within the FR \emph{in one step}, thus bypassing the remaining steps in the original procedure; 
2) An approach is proposed to select the superellipsoid ``squareness'' parameter to achieve a pre-specified level or near-optimality; 
3) Compared with the previous three-step approach, it is shown that a near-optimal solution, sufficiently close to \emph{global} optimality for practical purposes, can be guaranteed at the expense of a manageable increase in computation time. 

\color{black}
\section{RDOE via Superellipsoid-based Convex Optimisation}
\subsection{FR and the definition of superellipsoid}
As discussed in \cite{liu2022robust}, based on a UTOPF model with appropriate linearisation techniques, FR for active powers of all active customers can be compactly expressed as the following polyhedron. 
\begin{eqnarray}\small
\mathcal{F}(q)
    \label{fr-02}
    =\left\{p\left\vert\begin{matrix}
        Ap+Bq+Cv=d                 \\
        Ev\le f                    \\
    \end{matrix}\right.\right\}
	=\left\{p\left\vert\begin{matrix}
        \bar E Ap\le                  \\
        f+\bar E(d-Bq)                    \\
    \end{matrix}\right.\right\}
\end{eqnarray}
where 
$p$ and $q$ are variables corresponding to DOEs and reactive powers to be optimised, respectively;
$v$ is a vector consisting of variables related to nodal voltages;
$A,B,C,d,E$ and $f$ are constant parameters with appropriate dimensions 
and $\bar E=-EC^{-1}$.

The definition of the superellipsoid will be presented next, followed by the problem formulation in the next section. 

A $n$-dimension superellipsoid, with the centre being $u_c$ and lengths of all axes defined by a positive definite diagonal matrix $L$, can be expressed as \changed{$\mathcal{P}=\{p|p= Lw+u_c,||w||_n^n\le 1\}$ with $w_i$ being the $i^{th}$ element of $w$.} Fig.\ref{fig_superellipsoid} presents the boundary curves of $\mathcal{P}$ under various $n$, which shows that $\mathcal{P}$ is convex when $n\ge 1$ and becomes closer to the rectangle along with increasing value of $n$. 
\color{black}
Specially, when $n=2^K$ with $K$ being a positive integer, $\mathcal{P}$ can be expressed as $\mathcal{P}=\{p|p= Lw+u_c, ~w\in\mathcal{E}_h\}$, where $\mathcal{E}_h=\left\{y_1 |\begin{array}{cc}y_{k,i}^2\le y_{k+1,i} ~\forall i,\forall  k\le K-1;||y_K||_2\le 1\end{array}\right\}$ and $y_k$ is an intermediate vector variable with $y_{k,i}$ being its $i^{th}$ element. Moreover, $\mathcal{E}_h$ contains $mK+1$ quadratic constraints, where $m$ is the cardinality  of $y_k$.

\color{black}
\subsection{Mathematical formulation}
\color{black}
Noting that maximising the inscribed hyperrectangle regarding its \emph{volume}, say $V_h$, is equivalent to $\max_L\det(L)$ noting that 
\begin{enumerate}
    \item $V_h=2^v\prod\nolimits|L_{ii}w_i|=2^v\det(L)\prod\nolimits_i w_i$ with $v$ being the total number of active customers;
    \item $w_i^*=\arg\max\nolimits_{w:\norm{w}_n^n=1}{2^v\det(L)\sum\nolimits{\log(w_i)}}=v^{-1/n}$, the latter of which can be proved by applying the KKT condition when fixing $L$, leading to $V_h=(2v^{-1/n})^v\det(L)$.
\end{enumerate}

Then, $V_h$ can be maximised as follows by considering customers' operational statuses after replacing $n$ by $2^K$.
\color{black}
\begin{subequations}\small\label{max_hellip}
	\begin{eqnarray}
		\label{max_hellip_obj}
		\max_{\underline q\le q \le \bar q,L,u_c}{\log(\det (L))}-\varepsilon_{md}\sum\nolimits_i \delta_i\\
		s.t.~~
		\label{opt_01_rc}
            \bar EALy_1+\bar EAu_c\le f+\bar E(d-Bq)~~\forall y_1\in\mathcal{E}_h\\
  		\label{opt_02_rc}
		-\delta_i\le u_c(i)-\lambda_i L_{ii}v^{-1/2^K}\le \delta_i~~\forall i
	\end{eqnarray}
\end{subequations}
where
$\varepsilon_{md}$ is a penalty factor;
$\delta_i$ is a non-negative slack variable for customer $i$;
$u_c(i)$ is the $i^{th}$ element of vector $u_c$;
$\lambda_i$ is a binary number to indicate the operational status of customer $i$ with $\lambda_i=1/-1/0$ indicating \emph{import}/ \emph{export}/\emph{unknown} (could be \emph{import} or \emph{export}),
and $\underline q$ and $\bar q$ are the lower and upper limits for $q$, respectively. 

\color{black}
Obviously, as shown in \eqref{opt_02_rc}, the lower and upper capacity limits for customer $i$ are respectively $u_c(i)-L_{ii}v^{-1/2^K}$ and $u_c(i)+L_{ii}v^{-1/2^K}$, with $u_c(i)$ and $L_{ii}$ to be optimised.

\color{black}
Constraint \eqref{opt_02_rc}, together with \eqref{max_hellip_obj}, ensure the export/import limit is as close to 0 kW as possible when the customer is importing/exporting power so it can freely vary its power between 0 kW and the allocated capacity limit. For example, if customer $i$ is exporting power, $|u_c(i)+ {L_{ii}}v^{-1/2^K}|$, as its import limit, will be penalised in the objective function. Further, if a customer's operational status is unknown, $\sum\nolimits_i|u_c(i)|$ will be penalised in the objective function, leading to a DOE allocation where the import and export limits in their absolute values are as close to each other as possible. 

To derive the robust counterpart (RC) of \eqref{max_hellip} with respect to uncertain $y_1\in \mathcal{E}_h$ in \eqref{opt_01_rc}, the dual formulation of \eqref{max_hellip_rccons}, which is a generalised formulation for seeking the maximum value for the first term on the left-hand side of \eqref{opt_01_rc}, needs to be derived.
\begin{subequations}\small\label{max_hellip_rccons}
	\begin{eqnarray}
		\max\nolimits_{y_1,\cdots,y_K} x^Ty_1\\
		s.t.~~y_{k,i}^2\le y_{k+1,i} ~~\forall  k\le K-1,\forall i ~~(\alpha_{k,i})\\
		||y_K||_2\le 1~~(\alpha_K)
	\end{eqnarray}
\end{subequations}
where $x$ is a known vector representing a single row of $\bar EAL$, and $\alpha_{k,i}$ and $\alpha_K$ are Lagrangian multipliers.

Noting that the Lagrangian function for \eqref{max_hellip_rccons} is 
\begin{subequations}\footnotesize\label{max_hellip_rccons_L}
	\begin{eqnarray}
		W(\alpha_{k,i},\alpha_K,y_k)=-\sum_{2\le k\le K-1}\sum_i \left[\alpha_{k,i}y_{k,i}^2-\alpha_{k-1,i}y_{k,i}\right]\nonumber\\
		-\sum_i \left[\alpha_{1,i}y_{1,i}^2-x_iy_{1,i}\right]
		-\alpha_K(||y_K||_2-1)+\alpha_{K-1}^Ty_{K}\\
		\le \sum_i \frac{x_i^2}{4\alpha_{1,i}}+\sum_{2\le k\le K-1}\sum_i \frac{\alpha_{k-1,i}^2}{4\alpha_{k,i}}
		-\alpha_K(||y_K||_2-1)+\alpha_{K-1}^Ty_{K}
	\end{eqnarray}
\end{subequations}

Then $\min_{\alpha_{k,i},\alpha_K}\max_{y_k} W(\alpha_{k,i},\alpha_K,y_k)$, which is the equivalent problem of \eqref{max_hellip_rccons}, can be further formulated as 
\begin{subequations}\small\label{max_hellip_rccons_dual}
	\begin{eqnarray}
        \label{max_hellip_rccons_dual_01}
		\min_{\alpha_{k,i},\alpha_K} \alpha_K+\sum_i {x_i^2}/{4\alpha_{1,i}}+\sum_{2\le k\le K-1}\sum_i {\alpha_{k-1,i}^2}/{4\alpha_{k,i}}\\
		s.t.~~||\alpha_{K-1}||_2\le \alpha_K
	\end{eqnarray}
\end{subequations}

By introducing intermediate variables $t_{m,k,i}$, replacing $x$ by $[\bar EAL]_m^T$ and removing \emph{min} operator in \eqref{max_hellip_rccons_dual_01}, \eqref{max_hellip} can be reformulated as
\begin{subequations}\small\label{max_hellip_rc}
	\begin{eqnarray}
    \label{opt_final_obj}
		\max_{\underline q\le q \le \bar q,L,u_c}{\log(\det (L))}-\varepsilon_{md}\sum\nolimits_i \delta_i\\
		s.t.~\alpha_{m,K}+\sum\nolimits_{k\le K-1}\sum\nolimits_i t_{m,k,i}+[\bar EAu_c]_m \nonumber\\ \le [f+\bar Ed]_m-[\bar EBq)]_m~~\forall m\\
		([\bar EAL]_{m,i})^2\le 4\alpha_{m,1,i}t_{m,1,i}~\forall m,\forall i\\
		\alpha_{m,k-1,i}^2\le 4\alpha_{m,k,i}t_{m,k,i}~\forall m,\forall k\in\{2,\cdots,K-1\},\forall i\\
		||\alpha_{m,K-1}||_2\le \alpha_{m,K}~\forall m~\text{and}~\eqref{opt_02_rc}
	\end{eqnarray}
\end{subequations}
where $[\cdot]_m$ represents the $m^{th}$ row of a matrix or the $m^{th}$ element of a vector, and $[\cdot]_{m,i}$ represents the $(m,i)^{th}$ element of a matrix. 

Several remarks are given below.
\begin{enumerate}
    \item Since $\log(\det (L))$ is concave in $L$ and $x^2\le 4yz$ is equivalent to $ ||x;~y-z||_2\le y+z$, which is a second-order cone (SOC) constraint, \eqref{max_hellip_rc} is a convex optimisation problem. 
    \item For FR with $M$ inequalities and $v$ active customers, $M$ linear and $v\times M(K-1)+M$ quadratic constraints exist in \eqref{max_hellip_rc}, implying higher computational burden with larger $K$.
	\item \color{black}As the total DOE, defined as the absolute sum of all calculated capacity limits, is $S_h=2\sum\nolimits_i L_{ii} v^{-1/2^K}$ according to \eqref{opt_02_rc}, the maximum total DOE achievable is $\bar S_h=2\sum\nolimits_i L_{ii}$ when $K\rightarrow\infty$. With a given positive $\theta$ to make sure $1-S_h/\bar S_h\le\theta$, there is $K\ge(\log(\log v)-\log(-\log{(1-\theta)}))/\log{2}$. The value of $1-S_h/\bar S_h$ under various $K$ when $\theta=0.01$ is presented in Fig.\ref{fig_RDOE_best_K}, showing that setting $K$ as 10 can achieve acceptable accuracy with more than 500 active customers. 
    \begin{figure}[htb!]
    	\centering\includegraphics[scale=0.20]{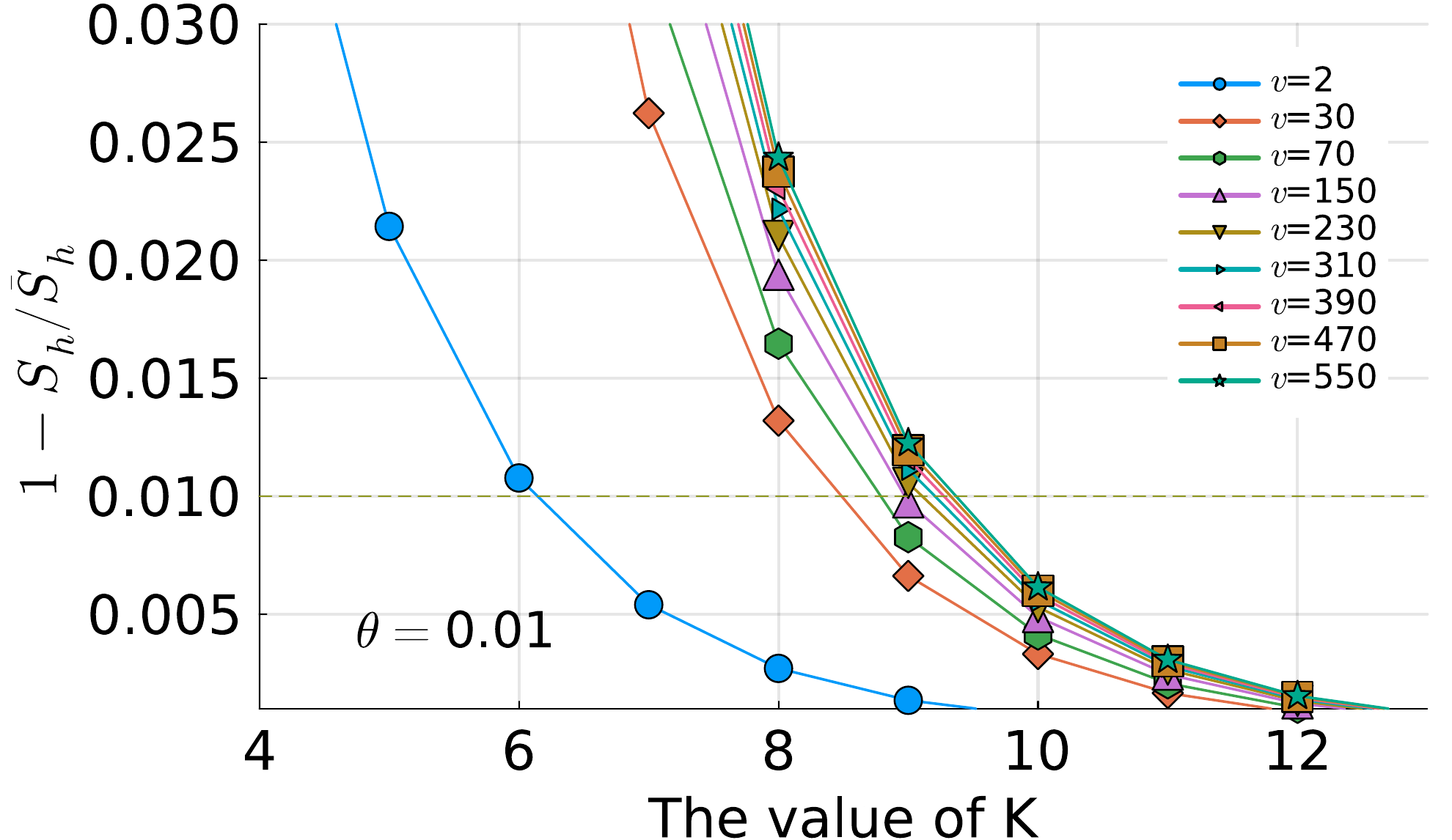}
    	\caption{\small The value of $1-S_h/\bar S_h$ with varying $K$.}
    	\label{fig_RDOE_best_K}
    \end{figure}


	\color{black}
    \item To further improve the computational efficiency, by introducing the intermediate variable $\gamma_{ii}$, the logarithmic term in \eqref{opt_final_obj} can be replaced by $\sum\nolimits_i \gamma_{ii}$ and an extra constraint $\gamma_{ii}\le ({L_{ii}-\bar L_{ii}^k})/{\bar L_{ii}^k}+\log\bar L_{ii}^k~(\forall k)$ with $\bar L_{ii}^k$ being the $k^{th}$ predefined point for the $i^{th}$ logarithmic term. 
    Generally, a larger number of predefined points, although leading to higher accuracy, results in higher computational complexity.
    \item \changed{The problem \eqref{max_hellip_rc} for both small-scale networks (the TwbAusNetworks in Section \ref{case_study}) and networks of larger scales (the AusNetwork and SynNetwork in Section \ref{case_study}) are solved by \texttt{Mosek} (version 10.0.40) on a desktop machine with Intel(R) Core(TM) i9-9900 CPU and 32 GB RAM after linearising \eqref{opt_final_obj}. For the three-step approach, \texttt{Ipopt} is employed for small-scale networks while \texttt{Mosek} and \texttt{Knitro} (trial version 13.2.0) are employed for others due to computational limitation of \texttt{Ipopt}, as discussed in \cite{liu2022robust}. Moreover, removing redundant constraints in the second step of the three-step approach is via \texttt{Xpress} (version 41.01.0).}
\end{enumerate}

\section{Case Study}\label{case_study}
\begin{figure}[htb!]
    \centering\includegraphics[scale=0.42]{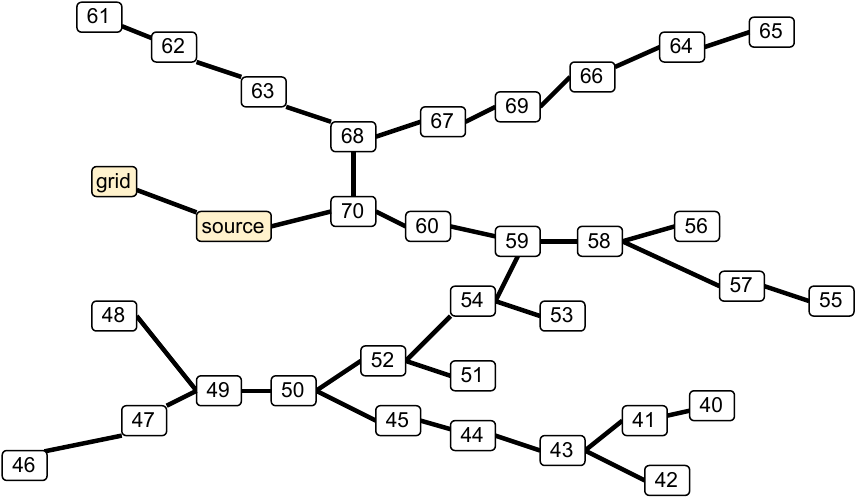}
    \caption{\small Topology of the AusNetwork (Information for customers' locations is available in \cite{LVFT_data} (\texttt{Network J}). Especially, customers 30, 21, 16, 29, 2, 18, 22, 19 and 17, which are constrained in Fig.\ref{fig_aus_J_hybrid_K_2_9_withQtrue_Mosek_SOCPtrue_LOBJtrue15}, are connected to buses 64, 40, 45, 65, 65, 41, 40, 46 and 42, respectively).}
    \label{fig_topology_simple_network_J}
\end{figure}

Three distribution networks are considered: a 2-bus illustrative
network (TwbNetwork), a 33-bus representative Australian network (AusNetwork) and a 132-bus synthetic network constructed by extending the AusNetwork (SynNetwork); see \cite{liu2022robust} for more details. For the 2-bus TwbNetwork, an ideal balanced voltage source with the voltage magnitude being $1.0~p.u.$ is connected to bus 1 and three customers (customers 1, 2 and 3) are connected to phase $b,a$ and $c$ of bus 2, respectively, where customer 2 is with fixed active/reactive power while DOEs are to be calculated for customers 1 and 3 with their reactive powers being controllable. The impedance of the three-phase line connecting bus 1 and bus 2 can be found in \cite{liu2022robust}. 
\changed{The AusNetwork/SynNetwork have 33/132 buses and 87/348 customers, of which 30/116 are active customers with controllable reactive powers. For the remaining passive customers, their active and reactive powers are fixed. The default limits on active and reactive powers for all studied networks are set as [-7 kW, 7 kW] and [-3 kvar, 3 kvar], respectively. The topology of the AusNetwork is presented in Fig.\ref{fig_topology_simple_network_J}, and the SynNetwork is generated by replicating the AusNetwork itself at buses ``40", ``55" and ``65". More details of the data, including network parameters, customers' locations and load profiles for passive customers, can be found in \cite{liu2022robust} and \cite{LVFT_data} (\texttt{Network J}) and are omitted here for simplicity.}

Simulation results for the TwbNetwork when customers 1 and 3 are importing powers are presented in Fig.\ref{fig_RDOE_concept_hypeellipse} and Table \ref{tab_case_study}. The deterministic approach (denoted as ``Dmtd") overestimates the DOEs compared with the globally optimal RDOEs (indicated by ``SO"). For the three-step approach (indicated by ``3step"), although RDOEs are ameliorated by the third step, noting the difference between the initial RDOEs (denoted by ``3step\_fnl") and the final RDOEs (indicated by ``3step\_ini"), it is still sub-optimal. By contrast, the superellipsoid-based approach (indicated by ``sESD") reports almost the same result as the SO-based approach when $\theta=0.01$ or $K=7$, demonstrating the efficiency of the proposed approach. 
\begin{figure}[htb!]
	\centering\includegraphics[scale=0.05]{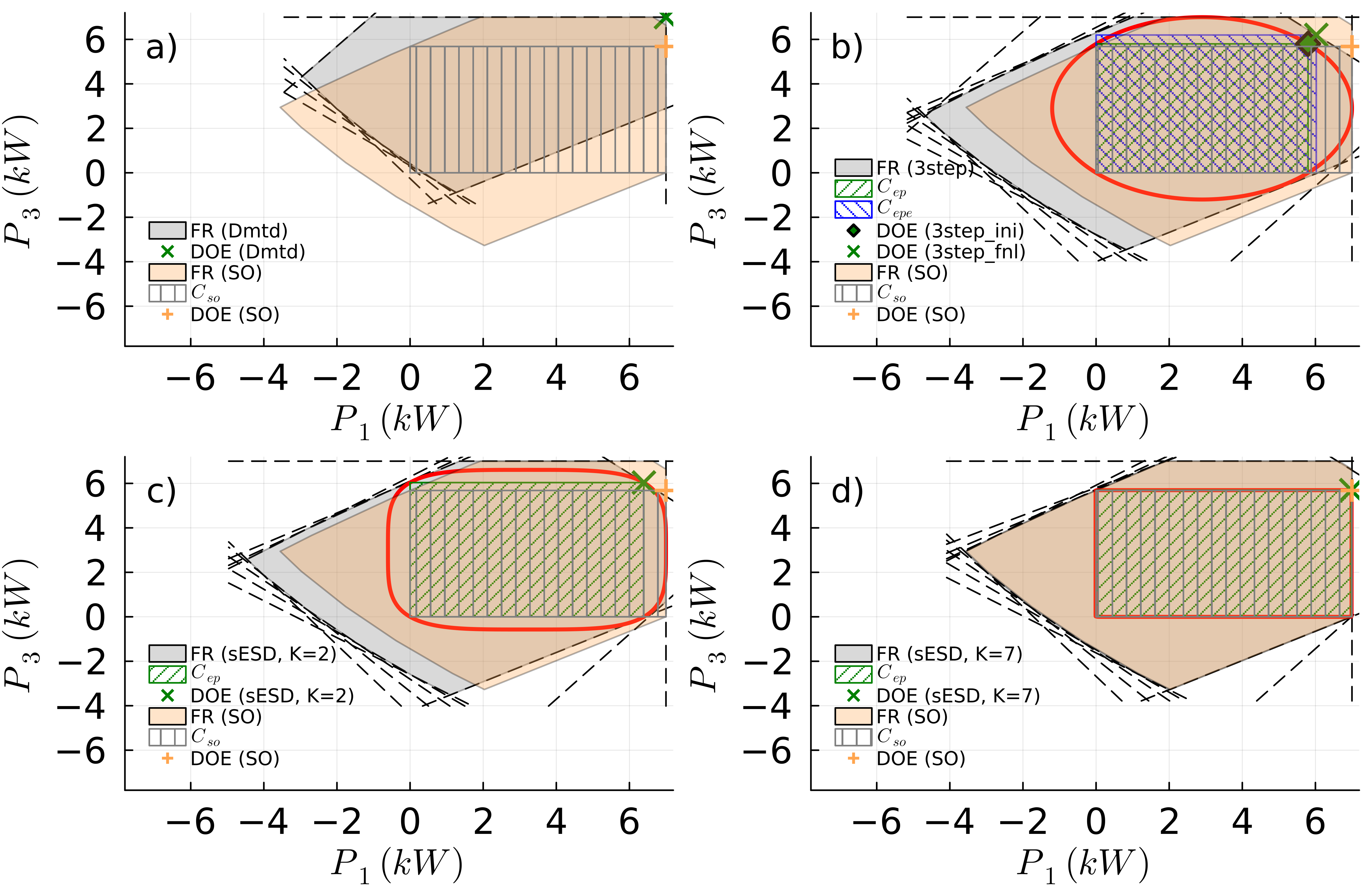}
	\caption{\small FRs and RDOEs compared with SO-based method under various approaches for TwbNetwork (boundaries for identified hyperellipsoid and superellipsoid are marked in red), where a): Dmtd; b): Three-step; c): sESD($K=2$); d): sESD($K=7$).}
	\label{fig_RDOE_2_bus_Qtrue_import_LOBJfalse}
\end{figure}

\begin{table}[htbp!]\footnotesize\renewcommand\arraystretch{1}
	\centering
	\setlength{\tabcolsep}{1.5pt}
	\caption{\small Total DOE calculated in kW under various approaches ($K^*$ set as 7/9/9 for TwbNetwork/AusNetwork/SynNetwork and the computational time in seconds are given under the calculated DOE in bracket).}
	\begin{tabular}{c|c|c|c|c|c}
		\hline\hline
		Network/Approach  & Dmtd & 3step & \makecell{sESD\\($K=2$)} & \makecell{sESD\\($K=K^*$)} & SO\\\hline
		TwbNetwork & \makecell{14.0\\(0.01 s)} & \makecell{12.2\\(0.10 s)} & \makecell{12.4\\(0.29 s)} & \makecell{12.7\\(0.94 s)} & \makecell{12.7\\(0.03 s)} \\\hline
		AusNetwork & \makecell{210.0\\(0.12 s)} & \makecell{149.5\\(28.46 s)} & \makecell{116.3\\(21.47 s)} & \makecell{165.8\\(102.14 s)} & --- \\\hline
		SynNetwork & \makecell{812.0\\(0.32 s)} & \makecell{230.8\\(303.39 s)}& \makecell{175.7\\(231.53 s)}& \makecell{269.0\\(1504.01 s)}&--- \\\hline		
		\hline
	\end{tabular}
	\label{tab_case_study}
\end{table}

\begin{figure}[htb!]
	\centering\includegraphics[scale=0.27]{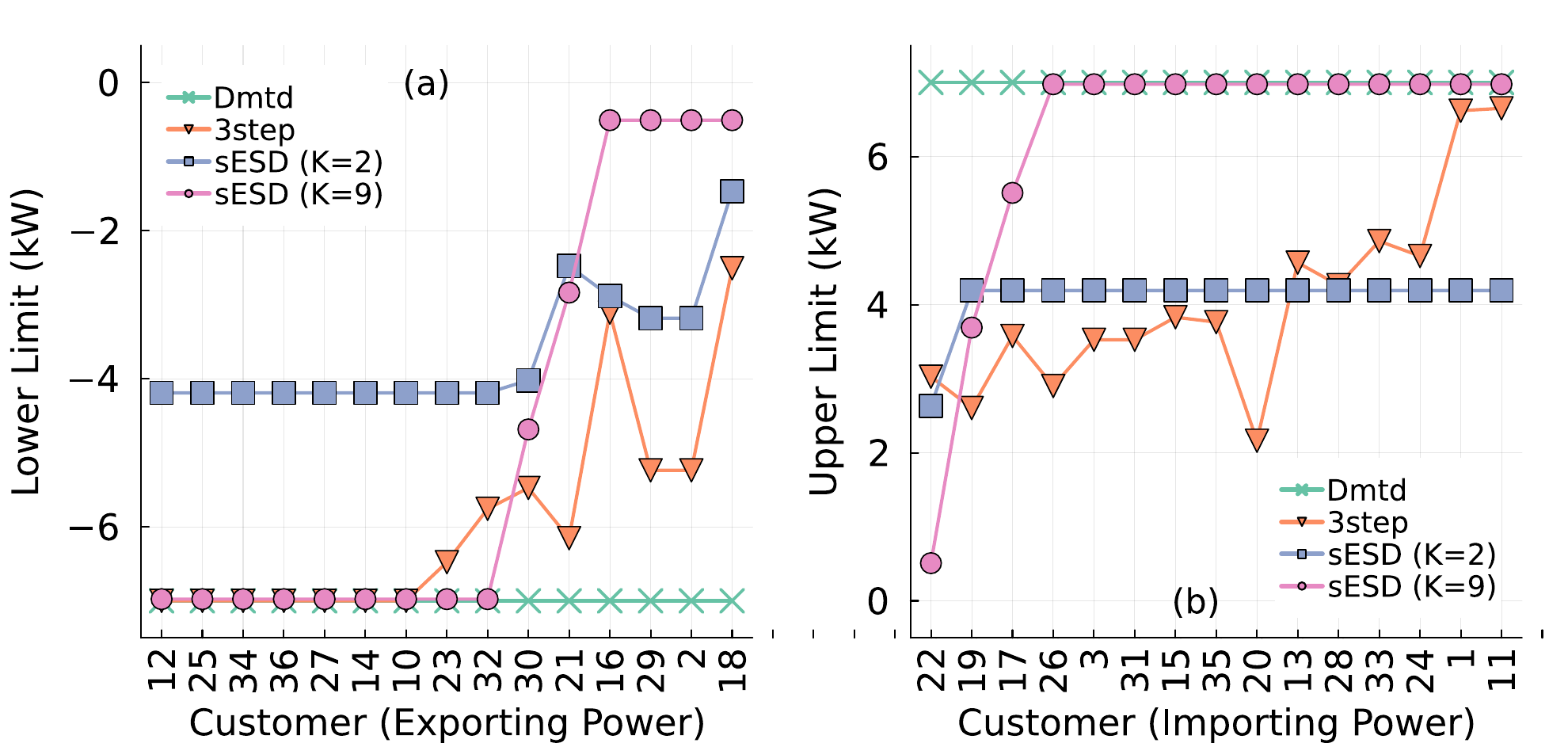}
	\caption{\small \changed{DOEs calculated by various approaches for the AusNetwork, where customers in sub-figure (a)/(b) are exporting/importing powers and their upper/lower limits are 0 kW.}}
	\label{fig_aus_J_hybrid_K_2_9_withQtrue_Mosek_SOCPtrue_LOBJtrue15}
\end{figure}

For the AusNetwork, the calculated DOEs are presented in Fig.\ref{fig_RDOE_2_bus_Qtrue_import_LOBJfalse} and Table \ref{tab_case_study}, where the deterministic approach again reports over-estimated results. \changed{Compared with the three-step approach, the total DOE calculated by the superellipsoid approach increases by 10.89\% to 165.8 kW when $K=9$ with, however, increased computational time, where the logarithmic terms are linearised with 15 predefined points. It is noteworthy that customers with locations at or near the network's \emph{endpoints} have smaller allocated capacity limits. This is because voltage levels at \emph{endpoints} are either higher or lower than those at \emph{internal} locations, making customers at or close to them more likely to be constrained.}
\begin{figure}[bht]
	\centering\includegraphics[scale=0.25]{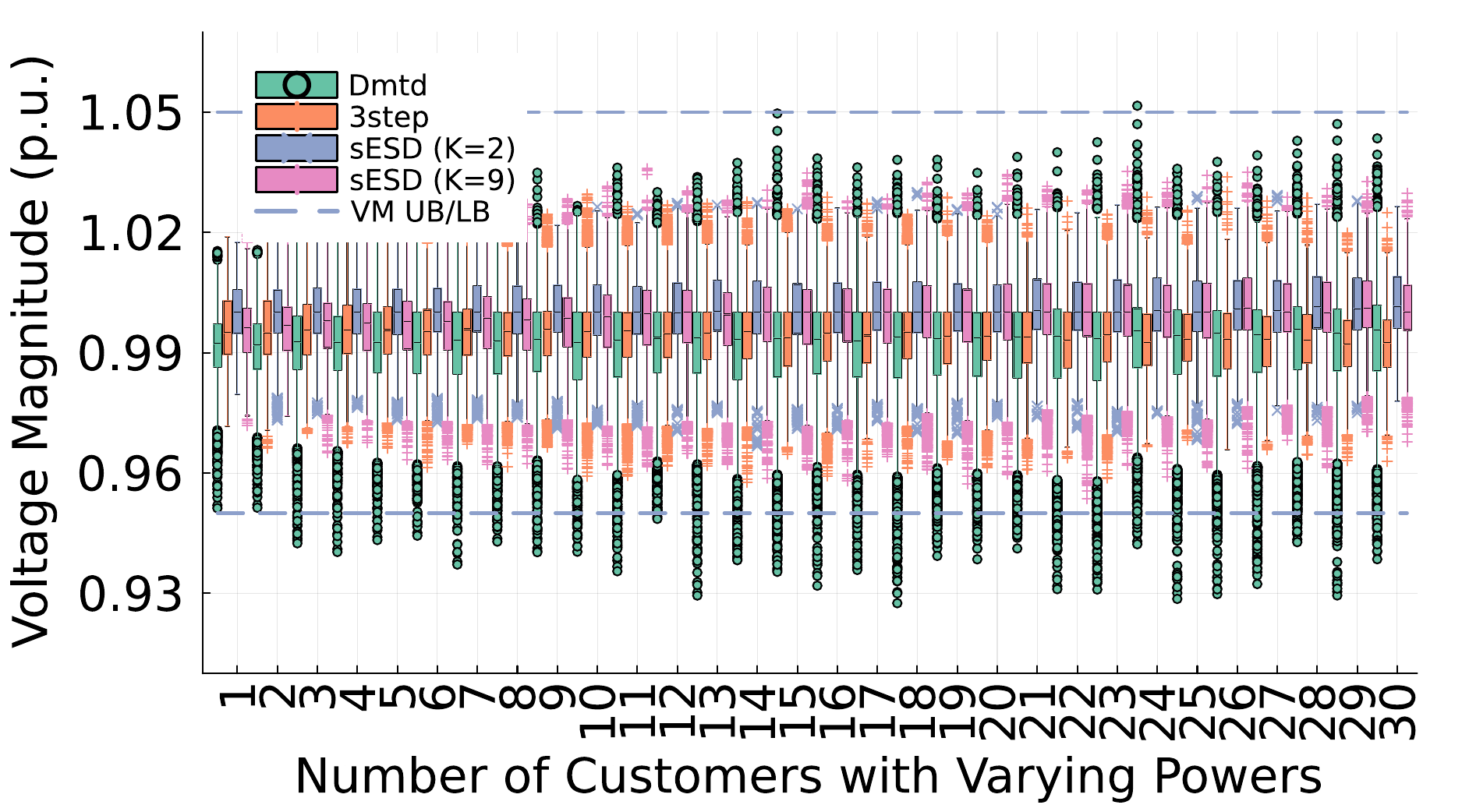}
	\caption{\small DOE assessment based on exact UTPF for the AusNetwork.}
	\label{fig_PFcompare_aus_J_hybrid_withQtrue_dftC7}
\end{figure}

Further assessment of the calculated DOEs by running exact UTPFs with
random load profiles\footnote{See \cite{liu2022robust} for details on generating the load profiles, here omitted for brevity.} on the AusNetwork is presented in Fig.\ref{fig_PFcompare_aus_J_hybrid_withQtrue_dftC7}, which clearly shows that the deterministic approach leads to significant voltage violations. In contrast, other approaches can guarantee operational security. 

\changed{The simulation results for the SynNetwork are presented in Table \ref{tab_case_study}, which again demonstrates that the deterministic approach overestimates the DOEs. In contrast, the superellipsoid-based approach increases the total DOE by 16.55\% compared with the three-step approach with $K=9$. Although the computational time increases to more than 1500 seconds, it is still acceptable when DOEs are calculated day-ahead, hourly or every several hours in intra-day operation. It is noteworthy that, since all extreme points of the hyperrectangle need to be enumerated in the SO approach, the constraints in the formulated problem would be more than $2^{30}>$1 \emph{billion} and $2^{116}$ for the AusNetwork and SynNetwork, respectively, which leads to computational intractability. Thus, the SO approach is not tested for the two networks.}

\section{Conclusions}

\changed{This letter proposes a superellipsoid-based approach to calculate RDOEs that can report a solution arbitrarily close to the global optimal via convex optimisation, in contrast to existing sub-optimal approaches. Increased computational time can be managed by selecting a proper parameter in formulating the superellipsoid while achieving the required accuracy, and in practice via a machine with higher computing power or a suitable time resolution, for example, every day, hourly or every several hours, to update the RDOEs if the computing power is limited. }

\color{black}

\bibliography{REFs_Power_Grid}

\end{document}